\documentclass[a4paper,12pt]{amsart}
\usepackage{fullpage}
\usepackage{amsmath}
\usepackage{amsfonts}
\usepackage{amssymb}

\usepackage[OT4]{fontenc}
\usepackage[latin2]{inputenc}
\DeclareMathOperator\dc{{\it i\partial\overline{\partial}}}
\DeclareMathOperator\co{{\it\mathbb C^n}}

\DeclareMathOperator\bk{{\it B_k}}
\DeclareMathOperator\om{{\it\omega}}
\DeclareMathOperator\Om{{\it\Omega}}

\DeclareMathOperator\Omb{{\it\Omega}^{''}}

\DeclareMathOperator \uk {{\it u_k}}
\DeclareMathOperator \ukk {{\it u_{k+1}}}
\DeclareMathOperator \wuk {{\it \widetilde{u}_k}}

\DeclareMathOperator \cuk {{\it \widehat{u}_k}}

\DeclareMathOperator \vk {{\it v_{k}}}
\DeclareMathOperator \omfr {{\it\om_{f}(d\rho^k,x_0)}}
\DeclareMathOperator \omfrkk {{\it\om_{f}(d\rho^{k+1},x_0)}}

\DeclareMathOperator \psh {{\it PSH(\Om)}}

\DeclareMathOperator \lbr {{\it\lbrace}}
\DeclareMathOperator \rbr {{\it\rbrace}}

\DeclareMathOperator\pa{\partial}

\def\om{\omega}
\def\al{\alpha}
\def\ga{\gamma}
\def\Om{\Omega}

\DeclareMathOperator \la {{\it\lambda}}

\newtheorem{theorem}{Theorem}

\newtheorem{remark}[theorem]{Remark}
\title{The $\mathcal C^{2,\alpha}$ estimate of complex Monge-Amp\`ere equation}
\subjclass[2000]{32W20}

\author{S\l awomir Dinew}
\address{S\l awomir Dinew,\\
Institute of Mathematics\\
Jagiellonian University, Poland}\email{slawomir.dinew@im.uj.edu.pl}
\author{Xi Zhang}
\address{Xi Zhang,\\
Department of Mathematics\\
Zhejiang University, P. R. China } \email{xizhang@zju.edu.cn}
\author{Xiangwen Zhang}
\address{Xiangwen Zhang,\\
Department of Mathematics and Statistics\\
McGill University, Canada} \email{xzhang@math.mcgill.ca}

\begin{document}
\maketitle
\begin{abstract}
 We prove that any $\mathcal C^{1,1}$ solution to complex Monge-Amp\`ere equa-\newline tion
 $det(u_{i\bar{j}})=f$ with
 $0<f\in\mathcal C^{\alpha}$ is in $\mathcal C^{2,\alpha}$ for $\alpha\in (0,1)$.
\end{abstract}
\section{introduction}
In a seminal paper \cite{Caf} Caffarelli proved, among other things, the following interior regularity for the Dirichlet problem for the {\it real} Monge-Amp\`ere equation:
\begin{theorem}\label{tw1} Le $\Om$ be a convex domain in $\mathbb R^n$ and u is a convex solution (understood in the viscosity sense) of the problem
\begin{equation}\label{caffa}
 det (u_{ij})=f,
\end{equation}
where $f$ is positive and $\alpha$-H\"older continuous for some $\alpha\in (0,1)$. Assume moreover that $u$ is equal to $0$ on $\pa\Om$. Then $u\in\mathcal C^{2,\alpha}(\Om)$.
\end{theorem}
The essence of this theorem is that under very mild assumptions on the right hand side we obtain in fact classical solutions.

It should be noted that, quite contrary to linear PDE problems there is no purely interior regularity for fully nonlinear PDE problems, even if we have analytic right hand side data. This follows from the famous Pogorelov example $$ u(x) = (x_1^2
+ 1)|x'|^{2\beta},$$
$$ det(u_{ij}) =c(n,\beta)(1 + x_1^2)^{n-2}((2\beta-1)-(2\beta+1)x_1^2)||x'||^{2(\beta n+1-n)},$$
where  $x=(x_1,x'),\ x'=(x_2,x_3,\cdots,x_n),\ \beta\geq1/2,\ n\geq 3$ (\cite{Po}). Thus assumptions on the smoothness of $\pa\Om$ as well as on $\phi$ are important parts of the data.

Recently Trudinger and Wang \cite{TW} have found the global version of Caffarelli's result. Their theorem reads as follows:
\begin{theorem}\label{tw2} Let $\Om$ be uniformly convex domain with $\mathcal C^3$ smooth boundary and $\phi\in\mathcal C^3(\bar{\Om})$. Let moreover $0<inf_{\Om}f\leq f\in \mathcal C^{\alpha}(\bar{\Om})$ Then any convex solution to the Dirichlet problem
\begin{equation*}
\begin{cases} 
  det (u_{ij})=f\ {\rm in}\ \Om,\\
u=\phi\ {\rm on}\ \pa\Om
 \end{cases}
\end{equation*}
satisfies the a priori estimate
$$||u||_{\mathcal C^{2,\alpha}(\bar{\Om})}\leq C,$$
where the constant $C$ depends on $n, \alpha, inf f, ||f||_{\mathcal C^{\alpha}(\bar{\Om})}, \pa\Om$ and $\phi$.

\end{theorem}
 
\begin{remark}
Both assumptions $\pa\Om\in \mathcal C^3$ and $\phi\in\mathcal C^3$ are sharp: see \cite{Wa2}.
\end{remark}

It should be also mentioned that if $f$ (still strictly positive) belongs to 
 $\mathcal C^{0,1}$, or even to some Sobolev space $W^{1,p}$ for $p>n$ 
  then one can use Evans-Krylov theory to the problem (see \cite{Ev,Kr}). Note however that generic
 H\"older function (except the Lipschitz ones) need not belong to any Sobolev space and thus the
 Evans-Krylov theory cannot be applied directly.

In this short note we investigate the interior regularity of the {\it complex} Monge-Amp\`ere equation with H\"older strictly positive right hand side data. Note that both Theorem \ref{tw1} and Theorem \ref{tw2} rely essentially on tools in convex analysis, like the geometric interpretation of the gradient image mappings, and {\it good shape} results for sublevel sets which are not available in the complex setting. Again due to lack of purely interior regularity one has to impose some additional assumptions on the solution $u$ itself. Here we shall work under the assumption that the solution is already in $\mathcal C^{1,1}(\Om)$. It remains an interesting problem to see where this condition can be weakened to $\Delta u\in L^{\infty}(\Om)$, or even more optimistically to $u\in\mathcal C^{1,\alpha}(\Om)$ for $\alpha>1-2/n$ (Pogorelov type examples in \cite{Bl} show that this exponent would be optimal).

Below we state our main result:
\begin{theorem}
 Let $\Om$ be a domain in $\co$ and $u\in\psh\cap\ \mathcal C^{1,1}(\Om)$ statisfy the Monge-Amp\`ere
equation
$$det(u_{i\bar{j}})=f.$$
Suppose additionally that $f\geq\la>0$ in $\Om$ for some constant $\la$ and $f\in\mathcal C^{\al}(\Om)$ for some $\al\in(0,1)$. Then $u\in\mathcal C^{2,\al}(\Om)$. Furthermore the $\mathcal C^{2,\al}$ norm of $u$ in any relatively compact subset is estimable in therms of $n, \al, \la, ||f||_{\mathcal C^{\al}(\Om)}$ and the distance of the set to $\pa\Om$.
\end{theorem}
 \begin{remark}
 In fact with a little more care the argument below (analogously to \cite{Wa1}) shows that if the right hand
 side is strictly positive and Dini continuous, then second derivatives of $u$ are Dini continuous too. Also
 if $f$ is strictly positive and Lipschitz then log-Lipshitz continuity of $D^2 u$ can be obtained.
\end{remark}

Quite recently Chen and He \cite{CH} proposed as an open problem to study the complex analogue of Theorem \ref{tw1} (that is the Dirichlet problem for the complex Monge-Amp\`ere equation with zero boundary data in a suitably smooth strictly pseudoconvex domain and strictly positive H\"older right hand side). We believe that our small contribution might be relevant to this issue.

An interesting element in our approach is that we have to work directly on the solution rather on its approximants.\footnote{For $u\in\mathcal C^{2}$ a more direct approach is possible by combining results below with a method from \cite{Lu}.} This at first sight causes some difficulties especially at the points where a priori $u$ does not have second order derivatives. It remains unclear to us if a sufficiently good approximation technique can be applied (one has to control simultaneously the approximants and their Monge-Amp\`ere functions).

In the proof we exploit a method due to Wang (\cite{Wa1}) (which originally comes from Caffarelli (\cite{Caf})) with small modifications from \cite{Ca}. We refer the reader also to the recent paper \cite{JW}, where Wang's approach is described in greater detail although the technical details are a bit different.\footnote{We wish to thank professor Xu-Jia Wang for pointing out this refence to us.} The crucial fact that is used instead of Wang's ''condition (A)'' (\cite{Wa1}) is the Bedford-Taylor interior $\mathcal C^{1,1}$ estimate (\cite{BT}). In fact this result is the sole reason why we need the assumption that $u\in\mathcal C^{1,1}$.

We also take the opportunity to describe Wang's method in the complex setting in greater detail, since its flexibility
 makes us believe that it can be applied in various different problems.

\section{preliminaries}
Below we collect several results that will be used in the note.

First of all recall that Rademacher theorem implies that $\mathcal C^{1,1}$ functions have almost everywhere defined second order derivatives. Moreover by classical distribution theory if $u\in \mathcal C^{1,1}$ then $\dc u$ is a current whose coefficient are $L^{\infty}$ functions coinciding with the (almost everywhere defined) mixed second order derivatives of $u$.

Recall below the comparison principle due to Bedford and Taylor:
\begin{theorem}[\cite{BT} Comparison principle]
Given a domain $\Om\subset\co$, let $u$ and $v$ be $\mathcal C^{1,1}(\Om)\cap \mathcal C(\bar{\Om})$ plurisubharmonic functions.\footnote{ Actually the theorem holds for merely locally bounded $u$ and $v$, see \cite{BT}. Here we state it in this form for the sake of simplicity.} Suppose that $u\leq v$ on $\pa \Om$ and $det(u_{i\bar{j}})\geq det(v_{i\bar{j}})$. Then $u\leq v$ in the whole $\Om$.
\end{theorem}

Building on this result and using the transitivity of the automorphism group of the unit ball in $\co$ Bedford and Taylor were able to prove the following interior a priori estimate:
\begin{theorem}[\cite{BT} Interior second order estimate]\label{c2}
 Let $\mathbb B$ be the unit ball in $\co$ and let $B'\subset\subset\mathbb B$ be arbitrary compact subset
 of $\mathbb B$. Let $u\in PSH(\mathbb B)\cap \mathcal C(\bar{\mathbb B})$ solve the Dirichlet problem
\begin{equation*}
 \begin{cases}
det(u_{i\bar{j}})=f,\\
u=\phi\ {\rm on}\ \pa\mathbb B,
 \end{cases}
\end{equation*}
where $\phi\in \mathcal C^{1,1}(\pa \mathbb B)$ and  $0\leq f^{1/n}\in C^{1,1}(\mathbb B)$. Then $u\in \mathcal C^{1,1}(\mathbb B)$ and moreover there exist a constant $C$\ dependent only on $dist\lbr B',\pa \mathbb B\rbr$ such that
$$||u||_{\mathcal C^{1,1}(B')}\leq C(||\phi||_{\mathcal C^{1,1}(\pa\mathbb B)}+ ||f^{1/n}||_{\mathcal C^{1,1}(\mathbb B)}).$$
\end{theorem}
\begin{remark}
 Note that no strict positivity of $f$ is needed. Observe also that this estimate is scale and translation
 invariant i.e. the same constant will work if we consider the Dirichlet problem in any ball with arbitrary
 radius (and suitably rescaled set $B'$).
\end{remark}
 
Finaly let us mention an interior $C^3$ estimate which (in the real case) is due to Calabi \cite{Cal} (the complex version due to Yau (\cite{Y}) for the global case and to Riebesehl and Schulz  (\cite{SR}) for a local estimate).

Here we state the complex version which will be the one we shall use:
\begin{theorem}[\cite{SR} Interior third order estimate]\label{c3}
Let $\Om$ be a domain in $\co$ and $u\in\psh\cap\ \mathcal C^{4}(\Om)$ statisfy the Monge-Amp\`ere
equation
$$det(u_{i\bar{j}})=f.$$
Then one has the interior third order estimate
$$||\nabla \Delta u||_{\Om'}\leq C$$
for a constant $C$ dependent only on $n, ||\nabla u||_{\Om}$, $||\Delta u||_{\Om}$, $\inf_{\Om}f$, $||D^1 f||_{\Om}$, $||D^2 f||_{\Om} $ and $dist\lbr \Om',\pa\Om\rbr$.
\end{theorem}

\medskip

\section{proof of main result}

\noindent{\bf proof of the main result:}
\par
 Fix $\Om''\Subset\Om'\Subset\Om$ and let $d:=dist(\Om'',\pa\Om')$. It is enough to show that $u\in\mathcal
 C^{2,\ga}(\Om'')$. As we have already mentioned we shall rely on the method developed by Wang (\cite{Wa1}).
 Some of the ideas are borrowed
 also from Campanato (\cite{Ca}).

For any fixed $x_0\in\Om''$ let us consider the system of balls $B(x_0,d\rho^k),\ k=0,\ 1,\ 2,\ \cdots$, where we put $\rho:=1/2$. Associated to any such system let $u(x;x_0,k)$ be the solution to the following Dirichlet problem
\begin{equation}\label{uk}
\begin{cases} u(x;x_0,k)\in{\it PSH}(B(x_0,d\rho^k))\cap\mathcal C(\bar{B}(x_0,d\rho^k)),\\
 det(u(x;x_0,k)_{i\bar{j}})=f(x_0)\ {\rm in}\ B(x_0,d\rho^k),\\
 u(x;x_0,k)=u\ {\rm on}\ \pa B(x_0,d\rho^k).
\end{cases}
\end{equation}
For notational ease we denote $B(x_0,d\rho^k)$ by $\bk$ and $u(x;x_0,k)$ by $\uk$. In case two different systems (with different centers) will appear, we shall mark them by $\wuk$ and $\cuk$ to make a distinction.

By the Bedford-Taylor interior $\mathcal C^{1,1}$ estimate applied to each ball $\bk$ we obtain that
$$||\uk||_{\mathcal C^{1,1}(B_{k+1})}\leq c_0(||u||_{\mathcal C^{1,1}(\pa B_{k})}+K),$$
where the constant $K$ depends merely on the supremum of $f$ on $\bk$. Since we assume that $u\in C^{1,1}(\Om)$ it easily follows that for some constant $c_1$ we have
\begin{equation}\label{c11}
||\uk||_{\mathcal C^{1,1}(B_{k+1})}\leq c_1||u||_{\mathcal C^{1,1}(\Om')}.
\end{equation}
Let us stress here that this estimate is independent of $\la$. Also $c_1$ does not depend on $k$ and $x_0$ since Bedford-Taylor estimate is scale invariant.

Now, since $f\geq \la$ in $\Om$ we obtain that on $\it{B}_{k+1}$
\begin{equation}\label{ellipticity}
 c_2(\la,||u||_{\mathcal C^{1,1}(\Om')})\dc||z||^2\leq\dc\uk\leq c_3(||u||_{\mathcal
 C^{1,1}(\Om')})\dc||z||^2.
\end{equation}

Observe that the argument above applies with no changes if instead of $u$ one uses its mollification $u^{(\epsilon)}$ for $\epsilon$ small enough, so that  everything is well defined in $\Om'$. As $||u^{(\epsilon)}||_{1,1}\rightarrow||u||_{1,1}$ we can temporarily work with $u^{(\epsilon)}$ (and we suppress the indice $\epsilon$ for the sake of readability). By the main theorem in \cite{CKNS} the solutions $\uk$ with the new boundary data coming from $u^{(\epsilon)}$ are smooth.
This allows one to apply the complex version of Calabi estimate to the problem (\ref{uk}). Thus,
  for any $\ga\in(0,1)$, we have
\begin{equation}\label{2gamma}
 ||\uk||_{\mathcal C^{2,\ga}(B_{k+2})}\leq c_4(c_2,c_3,d,n)/\rho^{k\ga}.
\end{equation}
Letting now $\epsilon\rightarrow 0^{+}$ we obtain that this estimate remains true for the original function $\uk$.

 Observe that the $(2,\gamma)$ norms of the functions $\uk$ may blow up at a controlled rate when
 $k\rightarrow\infty$.

Note that on $\it{B}_{k+3}$ the following holds
\begin{align*}
 0=log(f(x_0))-log(f(x_0))&=log(det(u_{k;i\bar{j}}))-log(det(u_{k+1;i\bar{j}}))=\\
\int_0^1\frac{d}{dt}log(det((t\uk&+(1-t)\ukk)_{i\bar{j}}))dt=\\
\sum_{i,j=1}^n\int_0^1(t\uk+(1-t&)\ukk)^{i\bar{j}}dt(\uk-\ukk)_{i\bar{j}},
\end{align*}
 where, as usual $(a^{i\bar{j}})$ denotes the inverse transposed matrix to $(a_{i\bar{j}})$, while
 $u_{k;i\bar{j}}$ is the $(i,\bar{j})$ th mixed complex derivative of $\uk$.

Thus we obtain that the difference $\vk:=\uk-\ukk$ satisfies on $\it{B}_{k+3}$ a linear elliptic equation
\begin{equation}\label{linear}
 \sum_{i,j=1}^nb^{i\bar{j}}\vk_{;{i\bar{j}}}=0.
\end{equation}

The coefficients $b^{i\bar{j}}$ satisfy, according to what we have proved so far, the estimates
$c_5||\zeta||^2\leq b^{i\bar{j}}\zeta_i\bar{\zeta}_{j}\leq c_6||\zeta||^2$ for any $\zeta\in\co$ and moreover $b^{i\bar{j}}$ are $\ga$-H\"older continuous with $||b^{i\bar{j}}||_{\mathcal C^{\ga}(\it{B}_{k+3})}\leq c_7/\rho^{k\ga}$.

This allows one to apply Schauder interior estimates to (\ref{linear}) (see \cite{GT}, Theorem 6.2). So
\begin{equation}\label{interior}
 ||\vk||_{\mathcal C^0(\it{B}_{k+3})}+sup_{\lbr x\in\it{B}_{k+3}\rbr}(d_{x}|D\vk(x)|)+sup_{\lbr
 x\in\it{B}_{k+3}\rbr}(d_{x}^2|D^2\vk(x)|)+
\end{equation}
\begin{align*}sup_{\lbr x,y\in\it{B}_{k+3}\rbr}(d_{x,y}^{2+\ga}&|D^2\vk(x)-D^2\vk(y)|/|x-y|^{\ga})\leq\\
c_8(c_5,c_6&,(c_7/\rho^{k\ga})(diam{\it B}_{k+3})^{\ga})||v_k||_{\mathcal C^0(\it{B}_{k+3})}.
\end{align*}
Here $d_z:=dist(z,\pa\it{B}_{k+3}),\ d_{z,w}=min\lbr d_z,\ d_w\rbr$, and $D\vk$ (respectively $D^2\vk$) denotes any first order (resp. second order) partial derivative of $\vk$.

We wish to point out that while $||b^{i\bar{j}}||_{\mathcal C^{\ga}(\it{B}_{k+3})}$ may blow up the term $(diam{\it B}_{k+3})^{\ga}$ compensates for this, so $c_8$ is an uniform constant independent of $k$.

The following estimates are straightforward consequences of (\ref{interior}):
\begin{equation}\label{6}
 ||\vk||_{\mathcal C^1(\it{B}_{k+4})}\leq (c_9/\rho^k)||v_k||_{\mathcal C^0(\it{B}_{k+3})};
\end{equation}
\begin{equation}\label{7}
 ||\vk||_{\mathcal C^2(\it{B}_{k+4})}\leq (c_{10}/\rho^{2k})||v_k||_{\mathcal C^0(\it{B}_{k+3})};
\end{equation}
\begin{equation}\label{8}
 ||\vk||_{\mathcal C^{2,\ga}(\it{B}_{k+4})}\leq (c_{11}/\rho^{(2+\ga)k})||v_k||_{\mathcal C^0(\it{B}_{k+3})}.
\end{equation}

Let now ${\it\om_{f}(r,x_0)}=osc_{B(x_0,r)}f$ denote the modulus of continuity of $f$ at $x_0$. Recall that $\uk$ and $u$ coincide on $\pa\bk$ and application of the comparison principle yields the inequalities
$$\uk(z)+\omfr(||z||^2-(d\rho^k)^2)\leq u(z)\leq \uk(z)-\omfr(||z||^2-(d\rho^k)^2)$$
for $z\in \bk$. Thus we get
\begin{equation}\label{9}
 ||\uk-u||_{\mathcal C^{0}(\bk)}\leq (d\rho^k)^2\omfr.
\end{equation}
Analogously
\begin{equation}\label{10}
 ||\ukk-u||_{\mathcal C^{0}({\it B}_{k+1})}\leq (d\rho^{k+1})^2\omfrkk,
\end{equation}
and coupling (\ref{9}) and (\ref{10}) we have
\begin{equation}\label{12}
 ||\vk||_{\mathcal C^{0}({\it B}_{k+4})}\leq ||\uk-u||_{\mathcal C^{0}({\it B}_{k+4})}+||\ukk-u||_{\mathcal
 C^{0}({\it B}_{k+4})}\leq c_{12}\rho^{2k}\omfr.
\end{equation}

Now, since we have assumed that $f\in \mathcal C^{\al}(\Om)$ we know that
\begin{equation}\label{alpha}
 \omfr\sim \rho^{\al}
\end{equation}
uniformly for $x_0\in\Om''$,
and thus (\ref{12}) together with (\ref{6}) and (\ref{7}) yield the fact that all the sequences $\lbr\uk(x_0)\rbr_k,\ \lbr D\uk(x_0)\rbr_k, \lbr D^2\uk(x_0)\rbr_k $ are Cauchy sequences and hence are all convergent. Analyzing the rate of convergence one easily sees that $\lim_{k\rightarrow\infty}\uk(x_0)=u(x_0)$ and $\lim_{k\rightarrow\infty}D\uk(x_0)=Du(x_0)$ for any first order partial derivative. The same in fact holds also for any second order partial derivative provided it exists at $x_0$. Since $u\in \mathcal C^{1,1}$ by assumption, it follows from Rademacher Theorem that this is the case almost everywhere. Note however that the limits $w_{st}(x):=\lim_{k\rightarrow\infty}\frac{\pa^2 u(x;x,k)}{\pa z_s\pa z_t}$, $\ w_{s\bar{t}}:=\lim_{k\rightarrow\infty}\frac{\pa^2 u(x;x,k)}{\pa z_s\pa \bar{z}_t}$ are defined in the whole $\Om''$.

Below we will show that all the limit functions $w_{st},\ w_{s\bar{t}}$ are in fact $\ga$-H\"older continuous. Thus second derivatives of $u$ exist almost everywhere and are equal to some $\ga$-H\"older continuous functions defined everywhere. By classical distribution teory it follows that $u\in \mathcal C^{2,\ga}$.

It is enough to prove H\"older continuity for any fixed $w_{st}$ (any $w_{s\bar{t}}$ goes the same way). To this end we fix two points $x,\ y\in \Omb$ and consider two cases:

{\it Case 1.} Let $||x-y||\geq d/16$. Then
$$|w_{st}(x)-w_{st}(y)|/||x-y||^{\al}\leq (16/d)^{\ga}|w_{st}(x)|+|w_{st}(y)|.$$
If $\wuk$ and $\cuk$ are the solutions of the Dirichlet problems related to systems of balls centered  at $x$ and $y$ respectively, by (\ref{7}) and (\ref{12}) we obtain that
\begin{equation}\label{wuk}
 |w_{st}(x)|=|\widetilde{u}_{0;st}(x)-\sum_{k=0}^{\infty}\widetilde{v}_{k;st}(x)|
\end{equation}
 with the obvious meaning of $\widetilde{v}_{k}$. The last quantity can be estimated as follows
$$|\widetilde{u}_{0;st}(x)-\sum_{k=0}^{\infty}\widetilde{v}_{k;st}(x)|\leq |\widetilde{u}_{0;st}(x)|+c_{14}\sum_{k=0}^{\infty}\rho^{k\al}\leq c_{15}<\infty,$$
where we have used (\ref{c11}) to control the first term.

Analogously on can bound $w_{st}(y)$ and thus in this case
\begin{equation}\label{16}
 |w_{st}(x)-w_{st}(y)|/||x-y||^{\al}\leq c_{16}.
\end{equation}
{\it Case 2.} Let now $||x-y||<d/16$. The we fix a $k\geq 0,\ k\in \mathbb N$ such that \newline $\rho^{k+5}d\leq||x-y||<\rho^{k+4}d$.

Then, as in \cite{Wa1} we estimate
\begin{align*}
|w_{st}(x)-w_{st}(y)|\leq& |w_{st}(x)-\widetilde{u}_{k;st}(x)|+\\
|w_{st}(y)-\widehat{u}_{k;st}(y)|+|\widehat{u}_{k;st}&(y)-\widetilde{u}_{k;st}(x)|=:I_1+I_2+I_3.
\end{align*}
The term $I_1$ can be easily handled in the following way:
\begin{equation*}
I_1=|\sum_{j=k}^{\infty}\widetilde{v}_{k;st}(x)|\leq \\
c_{17}\sum_{j=k}^{\infty}\rho^{j\al}\leq c_{18}(\rho^{k+5}d)^{\al}\leq c_{19}||x-y||^{\ga},
\end{equation*}
where we have used (\ref{7}), (\ref{12}) and (\ref{alpha}) in the first inequality.

$I_2$ is estimated completely the same way. To control $I_3$ observe that $y\in\it{B}(x,\rho^{k+4}d)$, so $\widetilde{u}_i$ is defined near $y$ for $i=0,\cdots, k$. Fix some $\ga>\al$ (the bigger the difference $\ga-\al$ the better).  We have

\begin{align*}
 I_3\leq
 |\widehat{u}_{k;st}(y)&-\widetilde{u}_{k;st}(y)|+|\widetilde{u}_{k;st}(y)-\widetilde{u}_{k;st}(x)|\leq\\
|\widehat{u}_{k;st}(y)-\widetilde{u}_{k;st}(y)|&+|\widetilde{u}_{0;st}(y)-\widetilde{u}_{0;st}(x)|+\sum_{j=0}^{k-1}|\widetilde{v}_{j;st}(y)-\widetilde{v}_{j;st}(x)|\leq\\
|\widehat{u}_{k;st}(y)-&\widetilde{u}_{k;st}(y)|+c_{20}|x-y|^{\gamma}\sum_{j=0}^{k-1}\rho^{(\al-\ga)j},
\end{align*}
where the last inequality holds because of (\ref{alpha}). Since $\al<\ga$ the last term is controlled by $||x-y||^{ga}\rho^{(\al-\ga)k}\sim||x-y||^{\al}$. Menawhile observe that $\wuk-\cuk$ satisfies in the domain $\it{B}(x,\rho^{k+2}d)\cap \it{B}(y,\rho^{k+2}d)$ a linear elliptic equation of type $c^{i\bar{j}}(\wuk-\cuk)_{i\bar{j}}=log(f(x))-log(f(y))$ analogous to the equation for $\vk$ (with $x$ and $y$ fixed here!). For the same reason as before we have uniform constants corresponding to $c_5,\ c_6,\ c_7$ (in fact if we have chosen at the beginning $c_5,\ c_6,\ c_7$ sufficiently big the same constants will do for this new equation). Thus Schauder interior estimates (with non homogeneous yet constant right hand side) give us the inequality
\begin{align*}
||\wuk-\cuk||_{\mathcal C^2(\it{B}(x,\rho^{k+3}d)\cap \it{B}(y,\rho^{k+3}d))}&\leq\\
c_{21}((1/\rho^{2k})||\wuk-\cuk||_{\mathcal C^0(\it{B}(x,\rho^{k+2}d)\cap\it{B}(y,\rho^{k+2}d))}&+|log(f(x))-log(f(y))|).
\end{align*}
Arguing as in (\ref{12}) and using that $f\geq\la>0$ the latter quantity is bounded by \newline $c_{22}\rho^{k\al}\leq c_{23}||x-y||^{\al}$, and that finishes the estimation of $I_3$.
Coupling all the obtained bounds we get the $\al$-H\"older continuity also in this case.
\medskip

{\bf Acknowledgements.} The note was finished while the first named author was visiting Princeton University. He wishes to thank this institution for the kind hospitality, and especially Professor Gang Tian for for his encouragement and help. Both the second and third named authors would like to thank Professor Pengfei Guan for the numerous helpful discussions on this problem and his constant encouragement and help. The note was written while the
 second named author was visiting McGill University. He would like to thank
 this institution for the hospitality.

%
%
%
%
%

\end{document}